\documentclass[11pt]{amsart}
\title{Towards commutator theory for relations}
\author{Paolo Lipparini}
\address{Dipartimento di Matematica, Viale della Ricca Scientifica,
II Universit\`a di Roma (Tor Vergata),
 ROME 
ITALY}
\thanks{The author has received support from MPI and GNSAGA. 
We thank H.-P. Gumm, K. A. Kearnes and E. W. Kiss 
for stimulating correspondence and discussions} 

\email{lipparin@axp.mat.uniroma2.it}
\urladdr{http://www.mat.uniroma2.it/\textasciitilde lipparin}
\newtheorem{Theorem}{Theorem}
\newtheorem{proposition}[Theorem]{Proposition}

\newtheorem{theorem}[Theorem]{Theorem}

\newtheorem{lemma}[Theorem]{Lemma}

\theoremstyle{definition}
\newtheorem{problem}[Theorem]{Problem}

\newcommand{\alg}{\mathbf} 

\begin{document}

\begin{abstract} 
In a general algebraic setting, 
we state some properties of commutators of 
reflexive admissible relations.
\end{abstract} 

\maketitle
\bigskip 

After commutator theory 
in Universal Algebra has been discovered 
about thirty years ago \cite{Sm}, 
many important results and applications
have been found. 
An introduction to commutator theory for
congruence modular varieties can be found in
 \cite{FMK} and \cite{Gu}. 
Shortly after, results valid for larger classes of varieties 
have been obtained in \cite{HMK} and \cite{LpTAMS}.
More recent results, as well as further references,
can be found, among others, in 
\cite{CHL,KS,LpCan2,LpActa}.

Present-day theory deals with commutators
of {\em congruences}. However, the possibility 
of a commutator theory for compatible reflexive relations
has been voiced
already in \cite[p. 186]{LpTAMS}.
Indeed, as noticed in 
\cite{Lpumi} 
(in part independently in \cite{Keordth}),
some notions
from classical commutator theory
can be extended to relations. 

If $\mathbf A$ is any algebra, and $ R, S$ are 
compatible and reflexive relations,
define $M( R , S )$ to be the set of all {\em matrices} of the form
\[
\begin{vmatrix}
t( \overline{a}, \overline{b})     & t( \overline {a}, \overline {b}\,') \cr
t( \overline {a}\,', \overline {b})     & t( \overline {a}\,', \overline {b}\,') 
\end{vmatrix}   
\] 
where $\overline{a}, \overline{a}\,' \in A^n $,
$\overline{b}, \overline{b}\,' \in A^m $, for some 
$m,n \geq 0$, $t$ is an $m+n$-ary term operation of $\mathbf A $,
and $\overline{a} R \overline{a}\,'$, $\overline{b} S \overline{b}\,'$.

We define $[ R, S ]$ to be the smallest
congruence that {\em centralizes} $R$ modulo $S$, that is,
the smallest congruence $ \delta $ such that $ z \delta w$ whenever
$x \delta y$ and  
$\begin{vmatrix}   
x & y \cr z & w 
\end{vmatrix}  
\in M( R , S ) 
$.

In our present setting, another commutator operation
is  more useful (cf. also \cite{Lpau}).
Let $[ R, S | 1]$
 be the {\em transitive closure} of the set   
\[\left\{ 
(z,w) | 
\begin{vmatrix}   
x & x \cr z & w 
\end{vmatrix}  
\in M( R , S ) 
\right\} 
\]

Notice that $[ R, S ]$ is, in general, much
larger than $[ R, S |1]$.
$[ R, S |1]$ is reflexive and compatible.
If $S$ is a tolerance, then  
$[ R, S | 1]$ is a congruence.
Clearly, $[ R, S | 1]$ is monotone in both arguments.
Moreover, $[ R, S | 1] \subseteq S^*$, and 
$[ R, S | 1] \subseteq Cg(R)$; actually,
$[ R, S | 1] \subseteq (S \cap (R^- \circ R))^*$.

For a relation $R$ on some algebra, let $R^\circ$ denote the smallest
tolerance containing $R$, and let $R^-$ denote the converse
of $R$. $R^*$ is the transitive closure of $R$, and
$Cg(R)$ is the smallest congruence containing $R$. 

Note the following easy but useful properties of $[ R, S | 1]$.

\begin{lemma}\label{x1a}
For $R, R_1,R_2,S,T,U$ reflexive and admissible relations on some algebra, the following hold: 

(i) $[R_1 \circ R_2,S|1] \subseteq 
\big(
S \cap ( R_2^- \circ (S \cap (R_1^- \circ R_1)) \circ R_2)
\big)^*$.

(ii) $ [R, S \circ T |1] \subseteq 
\big(
(R^- \circ R) \cap 
(
(S \cap (R^- \circ (S \cap T^-) \circ R))
\circ
(T \cap (R^- \circ (S^- \cap T) \circ R))
)
\big)^*$.

(iii) $ [R, S \circ T \circ U|1] \subseteq 
\big(
(R^- \circ R) \cap 
(S \circ (T \cap (R^- \circ (T \cap (S^- \circ U^-)) \circ R)) \circ U) 
\big)^*$. 
\end{lemma} 
 
\begin{proof}
Just draw a diagram.
\end{proof} 

Notice that it is possible to get a common generalization of (i) and (ii), as well as of
(i) and (iii). However, we get rather long formulae.

\begin{theorem}\label{x2} 
For every algebra $\alg A $, each of the following 
conditions imply all conditions below it:
\[ 
\tag{i}
R \subseteq [R, R^\circ|1]
\]
for every 
reflexive compatible relation $R$.
\[ 
\tag{ia}
R^* \subseteq [R, R^\circ|1]
\]
for every 
reflexive compatible relation $R$.
\[ 
\tag{ib}
R^-\subseteq [R, R^\circ|1]
\]
for every 
reflexive compatible relation $R$.
\[ 
\tag{ic}
R^\circ \subseteq [R, R^\circ|1]
\]
for every 
reflexive compatible relation $R$.
\[ 
\tag{id}
Cg(R) = [R, R^\circ|1]
\]
for every 
reflexive compatible relation $R$.
\[
\tag{ii}
R \cap T \subseteq [R, T |1]
\]
for every tolerance $T$ and every
reflexive compatible relation $R$.
\[
\tag{iii}
 (R_1 \circ R_2)\cap T \subseteq 
\big(
T \cap (R_2^- \circ (T \cap (R_1^- \circ R_1)) \circ R_2)
\big)^*
\]
for every tolerance $T$ and all
reflexive compatible relations $R_1$, $R_2$.
\[
\tag{iv}
 R_1 \cap (T \circ R_2) \subseteq
\big(
T \cap (R_2 \circ (T \cap (R_1^- \circ R_1)) \circ R_2^-)
\big)^*
\circ R_2
\] 
for every tolerance $T$ and all
reflexive compatible relations $R_1$, $R_2$.
\[
\tag{v}
\beta  \cap (T \circ S) \subseteq
\big(
T \cap (S \circ (T \cap \beta ) \circ S)
\big)^*
\circ S
\] 
for every congruence $ \beta $ and tolerances $T, S$.
\[
\tag{vi}
\beta  \cap (T \circ \gamma) \subseteq
\gamma \vee (T \cap \beta )^*
\] 
for every congruences $ \beta, \gamma $ and tolerance $T$.

Conditions
(i), (ia), (ib), (ic), (id) and (ii) above are equivalent for every algebra.
\end{theorem} 

\begin{proof}
(i) $ \Rightarrow $ (ia) 
$R \subseteq [R, R^\circ|1]$
implies $R^* \subseteq [R, R^\circ|1]^*= [R, R^\circ|1]$.

(ia) $ \Rightarrow $ (i) is trivial, since $R \subseteq R^*$.  

The proof that (i) and (ib)-(id) are equivalent is similar.

(i) $ \Rightarrow $ (ii)
Apply (i) with $ R \cap T $ in place of $R$, thus getting
\[
R \cap T \subseteq [R \cap T, (R \cap T)^\circ |1]\subseteq
[R , T^\circ|1] = [R , T|1]
\] 
since $[ , | 1]$ is monotone.

(ii) $ \Rightarrow $ (i)
By taking $T=R^\circ$ in (ii) we obtain (i).

(ii) $ \Rightarrow $ (iii) is immediate from Lemma
\ref{x1a}(i).

(iii) $ \Rightarrow $ (iv)
First notice that 
$R_1 \cap (T \circ R_2) \subseteq ((R_1 \circ R_2^-) \cap T) \circ R_2$,
since if $(a,c) \in R_1 \cap (T \circ R_2)$
then $a \, R_1 \, c$  and there is $b$ such that 
$ a\, T \, b \, R_2 \, c$, thus 
$ a\, R_1 \, c \, R_2^- \, b$ and 
$ (a,b) \in T \cap (R_1 \circ R_2^-)$. 

The conclusion follows by applying 
(iii) to $(R_1 \circ R_2^-) \cap T $. 

(iv) $ \Rightarrow $ (v) is trivial.

(v) $ \Rightarrow $ (vi) is trivial: take $S= \gamma$
and notice that
$\big(
T \cap (S \circ (T \cap \beta ) \circ S)
\big)^*
\subseteq 
\big(
S \circ (T \cap \beta ) \circ S
\big)^*=
\gamma \vee (T \cap \beta )^*
$
\end{proof}

By a similar argument, strengthening
condition (i) in Theorem \ref{x2}, we get:

\begin{theorem}\label{x3} 
For every algebra $\alg A $, each of the following 
conditions imply all conditions below it:
\[ 
\tag{i}
R \subseteq [R, R|1]
\]
for every 
reflexive compatible relation $R$.
\[ 
\tag{ia}
R^* \subseteq [R, R|1]
\]
for every 
reflexive compatible relation $R$.
\[
\tag{ii}
R \cap T \subseteq [R, T |1]
\]
for all reflexive compatible relations $T,R$.
\[
\tag{iii}
 (R_1 \circ R_2)\cap T \subseteq 
\big(
T \cap (R_2^- \circ (T \cap (R_1^- \circ R_1)) \circ R_2)
\big)^*
\]
for all
reflexive compatible relations $R_1$, $R_2$, $T$.
\[
\tag{iv}
 R_1 \cap (T \circ R_2) \subseteq
\big(
T \cap (R_2 \circ (T \cap (R_1^- \circ R_1)) \circ R_2^-)
\big)^*
\circ R_2
\] 
for all
reflexive compatible relations $R_1$, $R_2$, $T$.
\[
\tag{v}
\beta  \cap (T \circ S) \subseteq
\big(
T \cap (S \circ (T \cap \beta ) \circ S)
\big)^*
\circ S
\] 
for every congruence $ \beta $, tolerance $S$
and reflexive compatible relation $T$.
\[
\tag{vi}
\beta  \cap (T \circ \gamma) \subseteq
(\gamma \circ (T \cap \beta ))^*
\] 
for every congruences $ \beta, \gamma $
and reflexive compatible relation $T$.

Conditions
(i), (ia) and (ii) above are equivalent for every algebra.
\end{theorem} 

\begin{proposition}\label{x3a} 
For every algebra $\alg A $,  the following 
conditions are equivalent:
\[ 
\tag{i}
R \subseteq [R, R|1] \text{ and } R^- \subseteq [R, R|1]
\]
for every 
reflexive compatible relation $R$.
\[ 
\tag{ii}
Cg(R) = [R, R|1] 
\]
for every 
reflexive compatible relation $R$.
\end{proposition} 

There are other interesting consequences of
condition (i) in Theorem \ref{x2}, Theorem \ref{x3}, respectively. 
For example, we can 
apply conditions (ii) and (iii) in 
Lemma \ref{x1a}. As an example, we show:

\begin{theorem}\label{x4} 
(i) If an  algebra 
$\alg A $
satisfies
$R \subseteq [R, R^\circ|1]$
for every 
reflexive compatible relation $R$, 
then
$\alg A $
satisfies
\begin{multline*}
R \cap (S \circ T)
\cap (T^- \circ S^-)\subseteq
\\
\big(
(S \cap (R^- \circ (S \cap T^-) \circ R))
\circ
(T \cap (R^- \circ (S^- \cap T) \circ R))
\big)^*
\end{multline*}
for all 
reflexive compatible relations $R$, $S$, $T$.
In particular, if $\gamma$ is a congruence and
$ \gamma \supseteq S \cap T^-$ then
$\gamma \cap (S \circ T)
\cap (T^- \circ S^-)
\subseteq
((\gamma \cap S) \circ (\gamma \cap T))^*$.

(ii) If an  algebra 
$\alg A $
satisfies
$R \subseteq [R, R|1]$
for every 
reflexive compatible relation $R$, 
then
$\alg A $
satisfies
\begin{multline*}
R \cap (S \circ T)\subseteq
\big(
(S \cap (R^- \circ (S \cap T^-) \circ R))
\circ
(T \cap (R^- \circ (S^- \cap T) \circ R))
\big)^*
\end{multline*}
for all 
reflexive compatible relations $R$, $S$, $T$.
In particular, if $\gamma$ is a congruence and
$ \gamma \supseteq S \cap T^-$ then
$\gamma \cap (S \circ T)
\subseteq
((\gamma \cap S) \circ (\gamma \cap T))^*$.
\end{theorem} 

\begin{proof}
(ii) By the assumption with
$U=R \cap (S \circ T)$ in place of $R$, we have
$U
 \subseteq
 [U, U|1]
 \subseteq
 [R , S \circ T|1]$.
Now, apply
Lemma \ref{x1a}(ii).

The proof of (i) is similar.
Take
$U=R \cap (S \circ T)
\cap (T^- \circ S^-)
$ in place of $R$, thus getting
$U
 \subseteq
 [U, U^\circ|1]
 \subseteq
 [R , S \circ T|1]$,
since 
$U^\circ \subseteq
 ((S \circ T) \cap (T^- \circ S^-))^\circ=
 (S \circ T) \cap (T^- \circ S^-)\subseteq
 S \circ T$.
Again, apply
Lemma \ref{x1a}(ii).
\end{proof}

By using a more refined notation
(already introduced in \cite{Lpumi}),
 we can improve Lemma \ref{x1a}.

For $R, S, T$
compatible and reflexive relations,
let
\[
K( R, S; T) =
\left\{ 
(z,w) | 
\begin{vmatrix}   
x & y \cr z & w 
\end{vmatrix}  
\in M( R , S ), \ xTy 
\right\} 
\] 

Thus, 
 $[ R, S | 1]$
 is the transitive closure of 
 $K( R, S; 0)$, and  
 $[ R, S]$ 
is the smallest congruence
$ \gamma$ such that
 $K( R, S; \gamma) \leq \gamma$.  
Hence, the importance of the operator $K$
stems from the fact that any two elements
congruent modulo  $[ R, S]$ can be obtained
by a finite number of applications of 
 $K( R, S; - ) $ and of transitive closure 
and converse.  

\begin{lemma}
\label{x1b}
For $R, R_1,R_2,S,T,U, V$ reflexive and admissible relations on some 
algebra, the following hold: 

(i) 
$K(R_1 \circ R_2,S;V)
 \subseteq 
K( R_2,S; K(R_1,S;V))
$.

(ii) $K(R, S \circ T ;V) \subseteq 
K(R, S;S \cap (V \circ T^-)) 
\circ
K(R, T;T \cap (S^- \circ V)) 
$.

(iii) $K(R, S \circ T \circ U ;V) \subseteq 
(R^- \circ V \circ R) \cap 
(S \circ 
K(R,T; T \cap( S^-  \circ V \circ U^-)) 
\circ U) 
$. 
\end{lemma} 

Since, trivially, 
 $K( R, S; V) \subseteq  S \cap (R^- \circ (S \cap V) \circ R )$,
Lemma \ref{x1a} can be obtained as an immediate
consequence of Lemma 
\ref{x1b}, taking $V=0$.

\begin{problem}
Which of the following conditions
are equivalent {\em in a variety}?

(i) $R \subseteq [R, R^\circ|1]$
for every 
reflexive compatible relation $R$;

(ii) $R \subseteq [R, R|1]$
for every 
reflexive compatible relation $R$; 

(iii) $R \subseteq [R, R]$
for every 
reflexive compatible relation $R$;

(iv) $R \subseteq [R, R|1]$
for every 
tolerance $R$; 

(v) $R \subseteq [R, R]$
for every 
tolerance $R$.
\end{problem}

Notice that, for every algebra,
$R \subseteq [R, R]$
is equivalent to
$R \cap T \subseteq [R, T]$.
Announced results by K. Kearnes and E. Kiss suggest
the possibility that (i) and (v) above are not equivalent. 
Thus, probably, commutator theory for relations has stronger consequences
than commutator theory for tolerances, if we define the commutator
to be $[R, S|1]$ rather than $[R, S]$.

In a sequel to this paper we shall derive consequences from the 
existence of a difference term and of a weak difference term 
for $[R, R^\circ|1]$ and for
$[R, R|1]$. In particular, we shall deal with the following
properties

(a) $R \subseteq [R, R^\circ|1] \circ R^-$,

(b) $R \subseteq [R, R^\circ|1] \circ R^- \circ [R, R^\circ|1]$,

(c) $R \circ R \subseteq [R, R^\circ|1] \circ R$,

(d) $R \circ R \subseteq [R, R^\circ|1] \circ R \circ [R, R^\circ|1]$,

(e) $Cg(R)= [R, R^\circ|1] \circ R$,

(f) $Cg(R) = [R, R^\circ|1] \circ R \circ [R, R^\circ|1]$,
 
(g) $R \subseteq [R, R|1] \circ R^-$,

(h) $R \subseteq [R, R|1] \circ R^- \circ [R, R|1]$,

(i) $R \circ R \subseteq [R, R|1] \circ R$,

(j) $R \circ R \subseteq [R, R|1] \circ R \circ [R, R|1]$.

(k) $Cg(R)= [R, R|1] \circ R$,

(l) $Cg(R) = [R, R|1] \circ R \circ [R, R|1]$,

We shall also deal 
with a weaker commutator
\[
[ R, S | 1]_W=
\left\{ 
(x,w) | 
\begin{vmatrix}   
x & x \cr x & w 
\end{vmatrix}  
\in M( R , S ) 
\right\} 
\]

%
%
%
%


\def\cprime{$'$} \def\cprime{$'$}
\providecommand{\bysame}{\leavevmode\hbox to3em{\hrulefill}\thinspace}
\providecommand{\MR}{\relax\ifhmode\unskip\space\fi MR }
\providecommand{\MRhref}[2]{%
  \href{http://www.ams.org/mathscinet-getitem?mr=#1}{#2}
}
\providecommand{\href}[2]{#2}

\end{document}